\documentstyle[amscd,amssymb,verbatim,epsf]{amsart}

\begin{document}
\theoremstyle{plain}
\newtheorem{Thm}{Theorem}
\newtheorem{Cor}{Corollary}
\newtheorem{Con}{Conjecture}
\newtheorem{Main}{Main Theorem}
\newtheorem{Lem}{Lemma}
\newtheorem{Prop}{Proposition}

\theoremstyle{definition}
\newtheorem{Def}{Definition}
\newtheorem{Note}{Note}
\newtheorem{Ex}{Example}

\theoremstyle{remark}
\newtheorem{notation}{Notation}
\renewcommand{\thenotation}{}

\errorcontextlines=0
\numberwithin{equation}{section}
\renewcommand{\rm}{\normalshape}%

\title%
   {On the Space of Oriented Affine Lines in ${\Bbb{R}}^3$ }
\author{Brendan Guilfoyle}
\address{Brendan Guilfoyle\\
          Department of Mathematics and Computing \\
          Institute of Technology, Tralee \\
          Clash \\
          Tralee  \\
          Co. Kerry \\
          Ireland.}
\email{brendan.guilfoyle@@ittralee.ie}
\thanks{The first author was supported by the Isabel Holgate Fellowship from  Grey College, Durham and the Royal Irish Academy Travel Grant Scheme}
\author{Wilhelm Klingenberg}
\address{Wilhelm Klingenberg\\
 Department of Mathematical Sciences\\
 University of Durham\\
 Durham DH1 3LE\\
 United Kingdom.}
\email{wilhelm.klingenberg@@durham.ac.uk }

\date{December 7, 2002}

\begin{abstract}
We introduce a local coordinate description for the correspondence between 
the space of oriented affine lines in Euclidean ${\Bbb{R}}^3$ and the tangent bundle
to the 2-sphere. These can be utilised to give canonical coordinates on surfaces in ${\Bbb{R}}^3$, as we illustrate with a number of explicit examples.
\end{abstract}

\maketitle

The correspondence between oriented affine lines in ${\Bbb{R}}^3$ and the tangent bundle to the 2-sphere has a long history and has been used in various contexts. In particular,  it has been used in the construction of minimal surfaces \cite{weier}, solutions to the wave equation \cite{whitt} and the monopole equation \cite{hitch}. 

 The Euclidean group of rotations and translations acts upon the space of oriented lines ${\cal{L}}$ and in this paper we freeze out this group action by introducing a particular set of coordinates on ${\cal{L}}$. Our aim is to provide a local coordinate representation for the correspondence, thereby making it accessible to further applications.

One application is the construction of canonical coordinates on surfaces $S$ in ${\Bbb{R}}^3$ which come from the description of the normal lines of $S$ as local  sections of the tangent bundle of the 2-sphere. We illustrate this explicitly by considering the ellipsoid and the symmetric torus.

\begin{Def}
Let ${\cal{L}}$ be the set of oriented (affine) lines in Euclidean ${\Bbb{R}}^3$. 
\end{Def}
\begin{Def}
Let $\Phi:TS^2\rightarrow{\cal{L}}$ be the map that identifies ${\cal{L}}$ with the tangent bundle to
the unit 2-sphere in Euclidean ${\Bbb{R}}^3$, by parallel translation. This bijection gives ${\cal{L}}$ the structure of a differentiable 4-manifold.
\end{Def}

Let ($\xi,\eta$) be holomorphic coordinates on $TS^2$, where $\xi$ is obtained by 
stereographic projection from the south pole onto the plane through
the equator, and we identify ($\xi,\eta$) with the vector 
\[
\eta\frac{\partial}{\partial \xi}+\overline{\eta}\frac{\partial}{\partial \overline{\xi}}\in T_\xi S^2.
\]

\begin{Thm}
The map $\Phi$ takes ($\xi,\eta$)$\in TS^2$ to the oriented line given by
\begin{equation}\label{e:coord1}
z=\frac{2(\eta-\overline{\eta}\xi^2)+2\xi(1+\xi\overline{\xi})r}{(1+\xi\overline{\xi})^2}
\end{equation}
\begin{equation}\label{e:coord2}
t=\frac{-2(\eta\overline{\xi}+\overline{\eta}\xi)+(1-\xi^2\overline{\xi}^2)r}{(1+\xi\overline{\xi})^2},
\end{equation}
where $z=x^1+ix^2$, $t=x^3$, ($x^1,x^2,x^3$) are 
Euclidean coordinates on ${\Bbb{R}}^3={\Bbb{C}}\oplus{\Bbb{R}}$ and
$r$ is an affine parameter along the line such that $r=0$ is the
point on the line that lies closest to the origin.
\end{Thm} 
\begin{pf}
Stereographic projection from the south pole gives a map from ${\Bbb{C}}$ to ${\Bbb{R}}^3$ by
\begin{equation}\label{e:stereoproj}
z=\frac{2\xi}{1+\xi\overline{\xi}}\qquad\qquad\qquad
   t=\frac{1-\xi\overline{\xi}}{1+\xi\overline{\xi}}.
\end{equation}
The derivative of this map gives
\[
\frac{\partial}{\partial \xi}=\frac{2}{(1+\xi\overline{\xi})^2}\frac{\partial}{\partial z}-
   \frac{2\overline{\xi}\;^2}{(1+\xi\overline{\xi})^2}\frac{\partial}{\partial \overline{z}}-
   \frac{2\overline{\xi}}{(1+\xi\overline{\xi})^2}\frac{\partial}{\partial t},
\]
and similarly for its conjugate.

Thus
\begin{equation}\label{e:perpvec}
\eta\frac{\partial}{\partial \xi}+\overline{\eta}\frac{\partial}{\partial \overline{\xi}}=
\frac{2(\eta-\overline{\eta}\xi^2)}{(1+\xi\overline{\xi})^2}\frac{\partial}{\partial z} +
\frac{2(\overline{\eta}-\eta\overline{\xi}\;^2)}{(1+\xi\overline{\xi})^2}\frac{\partial}{\partial \overline{z}}-
\frac{2(\eta\overline{\xi}+\overline{\eta}\xi)}{(1+\xi\overline{\xi})^2}\frac{\partial}{\partial t}.
\end{equation}

Consider the line in ${\Bbb{R}}^3$ given by equations (\ref{e:coord1}) and (\ref{e:coord2}). The direction of this line is 
\[
\frac{2\xi}{1+\xi\overline{\xi}}\frac{\partial}{\partial z} +\frac{2\overline{\xi}}{1+\xi\overline{\xi}}\frac{\partial}{\partial \overline{z}}+\frac{1-\xi\overline{\xi}}{1+\xi\overline{\xi}}\frac{\partial}{\partial t}.
\]
When this unit vector is translated to the origin, it ends at the point $\xi\in S^2$ ({\it cf.} equation (\ref{e:stereoproj}))

The fixed vector determining the line is seen to be (\ref{e:perpvec}), and, using the fact that the Euclidean inner product of the basis vectors is
\[
\left(\frac{\partial}{\partial z},\frac{\partial}{\partial \overline{z}}\right)=\frac{1}{2} \qquad\qquad
\left(\frac{\partial}{\partial t},\frac{\partial}{\partial t}\right)=1
\]
\[
\left(\frac{\partial}{\partial z},\frac{\partial}{\partial z}\right)=
\left(\frac{\partial}{\partial \overline{z}},\frac{\partial}{\partial \overline{z}}\right)=
\left(\frac{\partial}{\partial z},\frac{\partial}{\partial t}\right)=
\left(\frac{\partial}{\partial \overline{z}},\frac{\partial}{\partial t}\right)=0,
\]
we compute that the line is orthogonal to the fixed
vector given by (\ref{e:perpvec}). Thus $r$ is an affine parameter
along the line such that $r=0$ is the point on the line that lies
closest to the origin, and the proof is completed. 
\end{pf}

Consider the map $\Psi:{\cal{L}}\times{\Bbb{R}}\rightarrow{\Bbb{R}}^3$ which takes a line and a number $r$ to a point on the line which is a parameter distance $r$ from the point on the line closest to the origin. 

\begin{Prop}
$\Psi^{-1}$ takes a point ($z,t$)$\in{\Bbb{R}}^3$ to a sphere in $\cal{L}\times{\Bbb{R}}$, the oriented lines containing the point: 
\[
\eta=\frac{1}{2}(z-2t\xi-\overline{z}\xi^2) \qquad\qquad r=\frac{\overline{\xi}z+\xi\overline{z}+(1-\xi\overline{\xi})t}{1+\xi\overline{\xi}}.
\]
\end{Prop}

\begin{pf}
This is comes from solving equations (\ref{e:coord1}) and
(\ref{e:coord2}) for $\eta$ and $r$. 

Alternatively, the second equation can be proved by finding the point
$p$ on the line with direction $\xi$ through ($z$, $t$)$\in R^3$ which minimises the distance to the origin. Then $r^2 = |(z,t)|^2 - |p|^2$, which gives the above expression for $r$. 

\end{pf}

By throwing away the $r$ information, the above formula gives the holomorphic sphere of lines through a given point ($z$,$t$)$\in{\Bbb{R}}^3$, as described in \cite{hitch}. These are a 3-parameter family of 
global sections of $TS^2$ and the associated line congruence in ${\Bbb{R}}^3$ is normal to round spheres about the given point.

More generally any oriented surface $S$ in ${\Bbb{R}}^3$ gives rise to a surface $\Sigma\subset\cal{L}$ through it's normal line congruence. Such a $\Sigma$ will, in general, not be holomorphic, nor be given by global sections of the bundle. However, locally, a surface can often be given by local non-holomorphic sections and the following examples illustrate this for two well-known surfaces.

The examples can be verified by substitution in equations (\ref{e:coord1}) and (\ref{e:coord2}) and then checking that the resulting surface, parameterised by its normal direction coordinate $\xi$, is indeed the one claimed.

\begin{Ex}

The triaxial ellipsoid with semi-axes $a_1$, $a_2$ and $a_3$ can be covered by coordinates $\xi$ via
\[
\eta=\frac{a_1(\xi+\overline{\xi})(1-\xi^2)+a_2(\xi-\overline{\xi})(1+\xi^2)
    -2a_3\xi(1-\xi\overline{\xi})}
     {2\sqrt{a_1(\xi+\overline{\xi})^2-a_2(\xi-\overline{\xi})^2
     +a_3(1-\xi\overline{\xi})^2}}
\]
\[
r=\sqrt{a_1\left(\frac{\xi+\overline{\xi}}{1+\xi\overline{\xi}}\right)^2
        -a_2\left(\frac{\xi-\overline{\xi}}{1+\xi\overline{\xi}}\right)^2
        +a_3\left(\frac{1-\xi\overline{\xi}}{1+\xi\overline{\xi}}\right)^2}.
\]
These coordinates extend to $\xi\rightarrow\infty$ and so this is an example of a global non-holomorphic section of $\pi:TS^2\rightarrow S^2$.
\end{Ex}

\begin{Ex}

The rotationally symmetric torus of radii $a$ and $b$ is given by
\[
\eta=\pm\frac{a}{2}\sqrt{\frac{\xi}{\overline{\xi}}}(1-\xi\overline{\xi})
\]
\[
r=b\pm\frac{2a\sqrt{\xi\overline{\xi}}}{1+\xi\overline{\xi}}.
\]

This describes the torus as a double cover of the 2-sphere, branched at the north and south poles.
\end{Ex}

\end{document}